\documentclass[12pt]{article}
\usepackage{amssymb}
\usepackage{amsmath}
\usepackage{amsfonts}
\usepackage[twocolumn=false]{geometry}

\input{tcilatex}
\begin{document}

\title{Subordination Properties of Certain Subclass of p-Valent Meromorphic
Functions Associated with Linear Operator}
\author{\textbf{R. M. El-Ashwah}$^{1}$\textbf{\ and A. H. Hassan}$^{2}$ \and 
$^{1}${\small Department of Mathematics, Faculty of Science, Damietta
University,} \and {\small New Damietta 34517, Egypt. E-mail:
relashwah@yahoo.com} \and $^{2}${\small Department of Mathematics, Faculty
of Science, Zagazig University,} \and {\small Zagazig 44519, Egypt.\ \
E-mail: alaahassan1986@yahoo.com}}
\date{}
\maketitle

\begin{abstract}
{\small In this paper, subordination results are studied for certain
subclass of p-valent meromorphic functions in the punctured unit disc having
a pole of order p at the origin. The subclass under investigation is defined
by using certain new linear operator. Moreover, we also introduced an
interesting particular cases of these results in several corollaries.}
\end{abstract}



\section*{1. Introduction}

Let $\Sigma _{p}$ denote the class of functions of the form%
\begin{equation}
f(z)=z^{-p}+\dsum\limits_{k=1-p}^{\infty }a_{k}z^{k}\text{ \ \ }(p\in 
%TCIMACRO{\U{2115} }%
%BeginExpansion
\mathbb{N}
%EndExpansion
:=\left\{ 1,2,3,...\right\} ),  \tag{1.1}
\end{equation}%
which are analytic in the punctured unit disc $U^{\ast }=U\backslash \left\{
0\right\} ;$ $U=\{z\in 
%TCIMACRO{\U{2102} }%
%BeginExpansion
\mathbb{C}
%EndExpansion
:|z|<1\}.$

\noindent For two functions $f(z)$ and $g(z),$ analytic in $U,$ we say that $%
f(z)$ is subordinate to $g(z)$ in $U,$ written $f\prec g$ or $f(z)\prec
g(z), $ if there exists a Schwarz function $\omega (z)$ which (by
definition) is analytic in $U,$ satisfying the following conditions (see
[12], [13]):%
\begin{equation*}
\omega (0)=0\text{ and }\left\vert \omega (z)\right\vert <1;\text{ }\left(
z\in U\right)
\end{equation*}%
such that 
\begin{equation*}
f(z)=g(\omega (z));\text{ }\left( z\in U\right) ,
\end{equation*}

\noindent Indeed it is known that%
\begin{equation*}
f(z)\prec g(z)\text{ \ \ }\left( z\in U\right) \Longrightarrow f(0)=g(0)%
\text{ \ \ and \ \ }f(U)\subset g(U).
\end{equation*}

\noindent In particular, If the function $g(z)$ is univalent in $U$, we have
the following equivalence (see also [4]):%
\begin{equation*}
f(z)\prec g(z)\text{ \ \ }\left( z\in U\right) \Longleftrightarrow f(0)=g(0)%
\text{ \ \ and \ \ }f(U)\subset g(U).
\end{equation*}

\noindent Following the recent work of El-Ashwah [7], for a function $f(z)$
in the class $\Sigma _{p}$, given by (1.1), the operator $L_{p}^{m}(\lambda
,\ell )$ is defined as following:%
\begin{equation}
L_{p}^{m}(\lambda ,\ell )f(z)=\left\{ 
\begin{array}{c}
f(z);\text{ \ \ \ \ \ \ \ \ \ \ \ \ \ \ \ \ \ \ \ \ \ \ \ \ \ \ \ \ \ \ \ \
\ \ \ \ \ \ \ \ \ \ }m=0 \\ 
\\ 
\tfrac{\ell }{\lambda }z^{-p-\tfrac{\ell }{\lambda }}\dint\limits_{0}^{z}t^{%
\left( \tfrac{\ell }{\lambda }+p-1\right) }L_{p}^{m-1}(\lambda ,\ell )f(t)dt;%
\text{ \ \ \ }m=1,2,...\text{ .}%
\end{array}%
\right.  \tag{1.2}
\end{equation}

\noindent Also, following the recent work of El-Ashwah and Hassan [9] (see
also [21]-[24]), for a function $f(z)\in \Sigma _{p}$, given by (1.1), also,
for $\mu >0,$ $a,c\in 
%TCIMACRO{\U{2102} }%
%BeginExpansion
\mathbb{C}
%EndExpansion
$ and $Re(c-a)\geq 0$, the integral operator%
\begin{equation*}
J_{p,\mu }^{a,c}:\Sigma _{p}\longrightarrow \Sigma _{p}
\end{equation*}%
is defined for $Re(c-a)>0$\ as follows:%
\begin{equation}
J_{p,\mu }^{a,c}f(z)=\frac{\Gamma (c-p\mu )}{\Gamma (a-p\mu )\Gamma (c-a)}%
\int\limits_{0}^{1}t^{a-1}(1-t)^{c-a-1}f(zt^{\mu })dt,  \tag{1.3}
\end{equation}%
and for $a=c$ as follows:%
\begin{equation}
J_{p,\mu }^{a,a}f(z)=f(z).  \tag{1.4}
\end{equation}

\noindent By iterations of the linear operators $L_{p}^{m}(\lambda ,\ell )$
defined by (1.2) and $J_{p,\mu }^{a,c}$ defined by (1.3) and (1.4), the
operator%
\begin{equation*}
I_{\lambda ,\ell }^{p,m}(a,c,\mu ):\Sigma _{p}\longrightarrow \Sigma _{p}
\end{equation*}%
is defined for the purpose of this paper by:%
\begin{equation}
I_{\lambda ,\ell }^{p,m}(a,c,\mu )f(z)=L_{p}^{m}(\lambda ,\ell )\left(
J_{p,\mu }^{a,c}f(z)\right) =J_{p,\mu }^{a,c}\left( L_{p}^{m}(\lambda ,\ell
)f(z)\right) .  \tag{1.5}
\end{equation}

\noindent Now, it is easily to see that the generalized operator $I_{\lambda
,\ell }^{p,m}(a,c,\mu )$ can be expressed as following:%
\begin{equation}
I_{\lambda ,\ell }^{p,m}(a,c,\mu )f(z)=z^{-p}+\frac{\Gamma (c-p\mu )}{\Gamma
(a-p\mu )}\dsum\limits_{k=1-p}^{\infty }\frac{\Gamma (a+\mu k)}{\Gamma
(c+\mu k)}\left[ \frac{\ell }{\ell +\lambda \left( k+p\right) }\right]
^{m}a_{k}z^{k},  \tag{1.6}
\end{equation}%
\begin{equation*}
\left( \mu {\small >}0;a,c{\small \in }%
%TCIMACRO{\U{2102} }%
%BeginExpansion
\mathbb{C}
%EndExpansion
,Re(a){\small >}p\mu ,Re(c{\small -}a){\small \geq }0;\ell {\small >}%
0;\lambda {\small >}0;m{\small \in }%
%TCIMACRO{\U{2115} }%
%BeginExpansion
\mathbb{N}
%EndExpansion
_{0}{\small =}%
%TCIMACRO{\U{2115} }%
%BeginExpansion
\mathbb{N}
%EndExpansion
\cup \left\{ 0\right\} ;p{\small \in }%
%TCIMACRO{\U{2115} }%
%BeginExpansion
\mathbb{N}
%EndExpansion
\right) .
\end{equation*}

\noindent In view of (1.2), (1.4) and (1.5), it is clear that:%
\begin{equation}
I_{\lambda ,\ell }^{p,0}(a,c,\mu )f(z)=J_{p,\mu }^{a,c}f(z)\text{ \ \ and \
\ }I_{\lambda ,\ell }^{p,m}(a,a,\mu )f(z)=L_{p}^{m}(\lambda ,\ell )f(z). 
\tag{1.7}
\end{equation}%
\noindent The operator $I_{\lambda ,\ell }^{p,m}(a,c,\mu )$ defined by (1.7)
has been extensively studied by many authors with suitable restrictions on
the parameters. For examples, see the following:

\noindent (i) $I_{\lambda ,\ell }^{1,-n}(a,c,\mu )=I_{\lambda ,\ell
}^{n}(a,c,\mu )f(z)\left( \mu {\small >}0;a,c{\small \in }%
%TCIMACRO{\U{2102} }%
%BeginExpansion
\mathbb{C}
%EndExpansion
,Re(c{\small -}a){\small \geq }0,Re(a){\small >}\mu ;\ell {\small >}%
0;\lambda {\small >}0;n{\small \in }%
%TCIMACRO{\U{2124} }%
%BeginExpansion
\mathbb{Z}
%EndExpansion
\right) $ (see El-Ashwah [6]);\smallskip

\noindent (ii) $I_{\lambda ,\ell }^{p,m}(p+\nu ,p+1,1)=I_{p,\nu
}^{m}(\lambda ,\ell )f(z)\left( m\in 
%TCIMACRO{\U{2115} }%
%BeginExpansion
\mathbb{N}
%EndExpansion
_{0};\lambda ,\ell ,\nu {\small >}0;p\in 
%TCIMACRO{\U{2115} }%
%BeginExpansion
\mathbb{N}
%EndExpansion
\right) $ (see El-Ashwah and Aouf [8]);\smallskip

\noindent (iii) $I_{\nu ,\lambda }^{1,m}(a+1,c+1,1)f(z)=\Im _{\lambda ,\nu
}^{m}(a,c)f(z)$ $\left( \lambda ,\nu {\small >}0;a\in 
%TCIMACRO{\U{2102} }%
%BeginExpansion
\mathbb{C}
%EndExpansion
;c\in 
%TCIMACRO{\U{2102} }%
%BeginExpansion
\mathbb{C}
%EndExpansion
\backslash 
%TCIMACRO{\U{2124} }%
%BeginExpansion
\mathbb{Z}
%EndExpansion
_{0}^{-};m\in 
%TCIMACRO{\U{2115} }%
%BeginExpansion
\mathbb{N}
%EndExpansion
_{0}\right) $ (see Raina and Sharma [16]);\smallskip

\noindent (iv) $I_{\lambda ,\ell }^{p,0}(a+p,c+p,1)f(z)=\ell _{p}(a,c)f(z)$ $%
\left( a\in 
%TCIMACRO{\U{211d} }%
%BeginExpansion
\mathbb{R}
%EndExpansion
;c\in 
%TCIMACRO{\U{211d} }%
%BeginExpansion
\mathbb{R}
%EndExpansion
\backslash 
%TCIMACRO{\U{2124} }%
%BeginExpansion
\mathbb{Z}
%EndExpansion
_{0}^{-},%
%TCIMACRO{\U{2124} }%
%BeginExpansion
\mathbb{Z}
%EndExpansion
_{0}^{-}=\left\{ 0,1,2,...\right\} ;p\in 
%TCIMACRO{\U{2115} }%
%BeginExpansion
\mathbb{N}
%EndExpansion
\right) $ (see Liu and Srivastava [11] and Srivastava and Patel
[18]);\smallskip

\noindent (v) $I_{1,\lambda }^{1,\beta }(\nu +1,2,1)f(z)=I_{\lambda ,\nu
}^{\beta }f(z)$ $\left( \beta \geq 0;\lambda >0;\nu >0\right) $ (see Piejko
and Sok\'{o}\l\ [15]);\smallskip

\noindent (vi) $I_{1,\lambda }^{1,n}(\nu +1,2,1)f(z)=I_{\lambda ,\nu
}^{n}f(z)$ $\left( n\in 
%TCIMACRO{\U{2115} }%
%BeginExpansion
\mathbb{N}
%EndExpansion
_{0};\lambda >0;\nu >0\right) $ (see Cho et al. [5]);\smallskip

\noindent (vii) $I_{\lambda ,\ell }^{1,0}(\nu +1,n+2,1)f(z)=\ell _{n,\nu
}f(z)$ $\left( n>-1;\nu >0\right) $ (see Yuan et al. [20]);\smallskip

\noindent (viii) $I_{\lambda ,\ell }^{p,0}(n+2p,p+1,1)f(z)=D^{n+p-1}f(z)$ $%
\left( n\text{ is an integer, }n>-p,\text{ }p\in 
%TCIMACRO{\U{2115} }%
%BeginExpansion
\mathbb{N}
%EndExpansion
\right) $ (see Uralegaddi and Somanatha [19], Aouf [1] and Aouf and
Srivastava [2]);\smallskip

\noindent (ix) $I_{1,1}^{p,\alpha }(a,a,\mu )f(z)=P_{p}^{\alpha }f(z)\left(
\alpha \geq 0;p\in 
%TCIMACRO{\U{2115} }%
%BeginExpansion
\mathbb{N}
%EndExpansion
\right) \ $(see$\ $Aqlan et al. [3]);\smallskip

\noindent (x) $I_{1,\beta }^{1,\alpha }(a,a,\mu )f(z)=P_{\beta }^{\alpha
}f(z)\left( \alpha ,\beta >0\right) $ (see Lashin [10]).

\section*{2. Preliminaries}

To establish our main results, we shall need the following lemmas:

\noindent \textbf{Lemma 1.} Using (1.6), we can obtain the following
recurrence relations of the operator $I_{\lambda ,\ell }^{p,m}(a,c,\mu )$:%
\begin{equation}
z\left( I_{\lambda ,\ell }^{p,m}(a,c,\mu )f(z)\right) ^{\prime }=\frac{%
a-p\mu }{\mu }I_{\lambda ,\ell }^{p,m}(a+1,c,\mu )f(z)-\frac{a}{\mu }%
I_{\lambda ,\ell }^{p,m}(a,c,\mu )f(z).  \tag{2.1}
\end{equation}%
and%
\begin{equation}
z\left( I_{\lambda ,\ell }^{p,m}(a,c+1,\mu )f(z)\right) ^{\prime }=\frac{%
c-p\mu }{\mu }I_{\lambda ,\ell }^{p,m}(a,c,\mu )f(z)-\frac{c}{\mu }%
I_{\lambda ,\ell }^{p,m}(a,c+1,\mu )f(z).  \tag{2.2}
\end{equation}%
Also,%
\begin{equation}
z\left( I_{\lambda ,\ell }^{p,m+1}(a,c,\mu )f(z)\right) ^{\prime }=\frac{%
\ell }{\lambda }I_{\lambda ,\ell }^{p,m}(a,c,\mu )f(z)-\frac{\ell +\lambda p%
}{\lambda }I_{\lambda ,\ell }^{p,m+1}(a,c,\mu )f(z).  \tag{2.3}
\end{equation}

\noindent \textbf{Lemma 2} [13]\textbf{.} Let the function $q(z)$ be
univalent in the unit disc $U$ and let $\theta $ and $\varphi $ are analytic
in a domain $D$ containing $q(U)$ with $q(w)\neq 0$ for all $w\in q(U)$. Set 
$Q(z)=zq^{\prime }(z)\varphi (q(z))$ and $h(z)=\theta (q(z))+Q(z)$. Suppose
that

\noindent (i) $Q(z)$ is starlike and univalent in $U;$

\noindent (ii) $\func{Re}\left\{ \frac{zh^{\prime }(z)}{Q(z)}\right\} >0$
for $z\in U$. If $p$ is analytic with $p(0)=q(0)$, $p(U)\subseteq D$ and%
\begin{equation}
\theta (p(z))+zp^{\prime }(z)\varphi (p(z))\prec \theta (q(z))+zq^{\prime
}(z)\varphi (q(z)),  \tag{2.4}
\end{equation}%
\noindent then%
\begin{equation}
p(z)\prec q(z)\text{ \ }(z\in U),  \tag{2.5}
\end{equation}%
\noindent and $q(z)$\ is\ the best dominant.\medskip

\noindent \textbf{Lemma 3} [17]\textbf{.} Let $q$ be a convex univalent
function in $U$ and let $\delta \in 
%TCIMACRO{\U{2102} }%
%BeginExpansion
\mathbb{C}
%EndExpansion
,$ $\gamma \in 
%TCIMACRO{\U{2102} }%
%BeginExpansion
\mathbb{C}
%EndExpansion
^{\ast }=%
%TCIMACRO{\U{2102} }%
%BeginExpansion
\mathbb{C}
%EndExpansion
\backslash \{0\}$ with%
\begin{equation}
\func{Re}\left\{ 1+\frac{zq^{\prime \prime }(z)}{q^{\prime }(z)}\right\}
>\max \left\{ 0,-\func{Re}\left\{ \frac{\delta }{\gamma }\right\} \right\} .
\tag{2.6}
\end{equation}%
\noindent If $p(z)$ is analytic in $U$ with $p(0)=q(0)$ and%
\begin{equation}
\delta p(z)+\gamma zp^{\prime }(z)\prec \delta q(z)+\gamma zq^{\prime }(z), 
\tag{2.7}
\end{equation}%
then%
\begin{equation}
p(z)\prec q(z)\text{ \ }(z\in U),  \tag{2.8}
\end{equation}%
\noindent and $q(z)$\ is\ the best dominant.\bigskip

\noindent In this paper, we find several sufficient conditions under which
some subordination results hold for the function $f\in \Sigma _{p}$ and for
suitable univalent function $q$ in $U$. We also introduced an interesting
particular cases of these results in several corollaries.

\section*{3. Subordination results}

Unless otherwise mentioned, we assume throughout the remainder of the paper
that $-1\leq B<A\leq 1$, $0\leq \alpha <p$, $\lambda >0$, $\ell >0$, $\mu >0$%
, $a,c\in 
%TCIMACRO{\U{2102} }%
%BeginExpansion
\mathbb{C}
%EndExpansion
$, $\func{Re}\left\{ a\right\} >p\mu $, $\func{Re}\left\{ c{\small -}%
a\right\} \geq 0$, $p\in 
%TCIMACRO{\U{2115} }%
%BeginExpansion
\mathbb{N}
%EndExpansion
$, $m\in 
%TCIMACRO{\U{2115} }%
%BeginExpansion
\mathbb{N}
%EndExpansion
_{0}$, $z\in U$ and the powers are principal.\medskip

\noindent We begin with investigating some sharp subordination results
regarding the operator $I_{\lambda ,\ell }^{p,m}(a,c,\mu )f(z).$

\noindent \textbf{Theorem 1. }Let $\xi \in 
%TCIMACRO{\U{2102} }%
%BeginExpansion
\mathbb{C}
%EndExpansion
^{\ast }=%
%TCIMACRO{\U{2102} }%
%BeginExpansion
\mathbb{C}
%EndExpansion
\backslash \left\{ 0\right\} $. Let the function $f\in \Sigma _{p}$ and the
function $q$ be univalent and convex in $U$ with $q(0)=1$. Suppose $f$ and $%
q $ satisfy any one of the following pairs of conditions:%
\begin{equation}
\func{Re}\left\{ 1+\frac{zq^{\prime \prime }(z)}{q^{\prime }(z)}\right\}
>\max \left\{ 0,-\frac{p}{\mu }\func{Re}\left\{ \frac{a-p\mu }{\xi }\right\}
\right\} ,\text{ \ \ \ \ \ \ \ \ \ \ \ \ \ \ \ \ \ \ \ \ \ \ \ \ \ \ \ \ \ \
\ \ \ \ \ \ }  \tag{3.1}
\end{equation}%
\begin{equation}
\frac{\xi }{p}\left( z^{p}I_{\lambda ,\ell }^{p,m}(a{\small +}1,c,\mu
)f(z)\right) {\small +}\frac{p{\small -}\xi }{p}\left( z^{p}I_{\lambda ,\ell
}^{p,m}(a,c,\mu )f(z)\right) \prec q(z){\small +}\frac{\mu \xi }{p\left( a%
{\small -}p\mu \right) }zq^{\prime }(z),  \tag{3.2}
\end{equation}

\noindent or

\begin{equation}
\func{Re}\left\{ 1+\frac{zq^{\prime \prime }(z)}{q^{\prime }(z)}\right\}
>\max \left\{ 0,-\frac{p}{\mu }\func{Re}\left\{ \frac{c-p\mu -1}{\xi }%
\right\} \right\} ,\text{ \ \ \ \ \ \ \ \ \ \ \ \ \ \ \ \ \ \ \ \ \ \ \ \ \
\ \ \ \ \ \ \ \ \ \ \ \ \ }  \tag{3.3}
\end{equation}%
\begin{equation}
\frac{\xi }{p}\left( z^{p}I_{\lambda ,\ell }^{p,m}(a,c{\small -}1,\mu
)f(z)\right) {\small +}\frac{p{\small -}\xi }{p}\left( z^{p}I_{\lambda ,\ell
}^{p,m}(a,c,\mu )f(z)\right) \prec q(z){\small +}\frac{\mu \xi }{p\left( c%
{\small -}p\mu {\small -}1\right) }zq^{\prime }(z),  \tag{3.4}
\end{equation}%
\noindent or%
\begin{equation}
\func{Re}\left\{ 1+\frac{zq^{\prime \prime }(z)}{q^{\prime }(z)}\right\}
>\max \left\{ 0,-\frac{p\ell }{\lambda }\func{Re}\left\{ \frac{1}{\xi }%
\right\} \right\} ,\text{ \ \ \ \ \ \ \ \ \ \ \ \ \ \ \ \ \ \ \ \ \ \ \ \ \
\ \ \ \ \ \ \ \ \ \ \ \ \ \ \ \ \ \ \ }  \tag{3.5}
\end{equation}%
\begin{equation}
\frac{\xi }{p}\left( z^{p}I_{\lambda ,\ell }^{p,m-1}(a,c,\mu )f(z)\right) 
{\small +}\frac{p{\small -}\xi }{p}\left( z^{p}I_{\lambda ,\ell
}^{p,m}(a,c,\mu )f(z)\right) \prec q(z){\small +}\frac{\lambda \xi }{\ell p}%
zq^{\prime }(z)\text{. \ \ \ \ \ \ \ \ \ }  \tag{3.6}
\end{equation}

\noindent Then%
\begin{equation}
z^{p}I_{\lambda ,\ell }^{p,m}(a,c,\mu )f(z)\prec q(z),  \tag{3.7}
\end{equation}

\noindent and $q(z)$ is the best dominant of (3.7).

\noindent \textbf{Proof.} Let%
\begin{equation}
k(z)=z^{p}I_{\lambda ,\ell }^{p,m}(a,c,\mu )f(z),  \tag{3.8}
\end{equation}

\noindent then it is easily to show that $k(z)$ is analytic in $U$ and $%
k(0)=1$. Differentiating both sides of (3.8) with respect to $z$, followed
by applications of the identities (2.1), (2.2) and (2.3) yield respectively%
\begin{equation}
z^{p}I_{\lambda ,\ell }^{p,m}(a+1,c,\mu )f(z)=k(z)+\frac{\mu }{a-p\mu }%
zk^{\prime }(z),\text{ \ \ \ \ \ }  \tag{3.9}
\end{equation}%
\begin{equation}
z^{p}I_{\lambda ,\ell }^{p,m}(a,c-1,\mu )f(z)=k(z)+\frac{\mu }{c-p\mu -1}%
zk^{\prime }(z),\text{ }  \tag{3.10}
\end{equation}

\noindent and%
\begin{equation}
z^{p}I_{\lambda ,\ell }^{p,m-1}(a,c,\mu )f(z)=k(z)+\frac{\lambda }{\ell }%
zk^{\prime }(z)\text{. \ \ \ \ \ \ \ \ \ \ \ \ \ \ }  \tag{3.11}
\end{equation}

\noindent Now, the subordination conditions (3.2), (3.4), and (3.6) are
respectively equivalent to%
\begin{equation}
k(z)+\frac{\mu \xi }{p\left( a{\small -}p\mu \right) }zk^{\prime }(z)\prec
q(z)+\frac{\mu \xi }{p\left( a{\small -}p\mu \right) }zq^{\prime }(z),\text{
\ \ \ \ \ \ \ \ }  \tag{3.12}
\end{equation}%
\begin{equation}
k(z)+\frac{\mu \xi }{p\left( c{\small -}p\mu -1\right) }zk^{\prime }(z)\prec
q(z)+\frac{\mu \xi }{p\left( c{\small -}p\mu -1\right) }zq^{\prime }(z), 
\tag{3.13}
\end{equation}

\noindent and%
\begin{equation}
k(z)+\frac{\xi \lambda }{p\ell }zk^{\prime }(z)\prec q(z)+\frac{\xi \lambda 
}{p\ell }zq^{\prime }(z)\text{. \ \ \ \ \ \ \ \ \ \ \ \ \ \ \ \ \ \ \ \ \ \
\ \ \ }  \tag{3.14}
\end{equation}

\noindent Therefore, applying Lemma 3 to each of the subordination
conditions (3.12), (3.13) and (3.14) with appropriate choices of $\delta $
and $\gamma $ we get the assertion (3.7) of Theorem 1. Then the proof of
Theorem 1 is completed.\bigskip

\noindent Putting $q(z)=\tfrac{1+Az}{1+Bz}$ in Theorem 1, we obtain the
following corollary:

\noindent \textbf{Corollary 1.} Let $\xi \in 
%TCIMACRO{\U{2102} }%
%BeginExpansion
\mathbb{C}
%EndExpansion
^{\ast }$. Let the function $f\in \Sigma _{p}$. Suppose any one of the
following pairs of conditions is satisfied:%
\begin{equation}
\tfrac{\left\vert B\right\vert -1}{\left\vert B\right\vert +1}<\frac{p}{\mu }%
\func{Re}\left\{ \tfrac{a-p\mu }{\xi }\right\} ,\text{ \ \ \ \ \ \ \ \ \ \ \
\ \ \ \ \ \ \ \ \ \ \ \ \ \ \ \ \ \ \ \ \ \ \ \ \ \ \ \ \ \ }  \tag{3.15}
\end{equation}%
\begin{equation}
\tfrac{\xi }{p}\left( z^{p}I_{\lambda ,\ell }^{p,m}(a{\small +}1,c,\mu
)f(z)\right) {\small +}\tfrac{p{\small -}\xi }{p}\left( z^{p}I_{\lambda
,\ell }^{p,m}(a,c,\mu )f(z)\right) \prec \tfrac{1+Az}{1+Bz}{\small +}\tfrac{%
\mu \xi }{p\left( a{\small -}p\mu \right) }\tfrac{\left( A{\small -}B\right)
z}{\left( 1{\small +}Bz\right) ^{2}},  \tag{3.16}
\end{equation}

\noindent or

\begin{equation}
\tfrac{\left\vert B\right\vert -1}{\left\vert B\right\vert +1}<\frac{p}{\mu }%
\func{Re}\left\{ \tfrac{c-p\mu -1}{\xi }\right\} ,\text{ \ \ \ \ \ \ \ \ \ \
\ \ \ \ \ \ \ \ \ \ \ \ \ \ \ \ \ \ \ \ \ \ \ \ \ \ \ \ \ }  \tag{3.17}
\end{equation}%
\begin{equation}
\tfrac{\xi }{p}\left( z^{p}I_{\lambda ,\ell }^{p,m}(a,c{\small -}1,\mu
)f(z)\right) {\small +}\tfrac{p{\small -}\xi }{p}\left( z^{p}I_{\lambda
,\ell }^{p,m}(a,c,\mu )f(z)\right) \prec \tfrac{1+Az}{1+Bz}{\small +}\tfrac{%
\mu \xi }{p\left( c{\small -}p\mu {\small -}1\right) }\tfrac{\left( A{\small %
-}B\right) z}{\left( 1{\small +}Bz\right) ^{2}},  \tag{3.18}
\end{equation}%
\noindent or%
\begin{equation}
\tfrac{\left\vert B\right\vert -1}{\left\vert B\right\vert +1}<\frac{p\ell }{%
\lambda }\func{Re}\left\{ \frac{1}{\xi }\right\} ,\text{ \ \ \ \ \ \ \ \ \ \
\ \ \ \ \ \ \ \ \ \ \ \ \ \ \ \ \ \ \ \ \ \ \ \ \ \ \ \ \ \ \ \ \ \ } 
\tag{3.19}
\end{equation}%
\begin{equation}
\tfrac{\xi }{p}\left( z^{p}I_{\lambda ,\ell }^{p,m-1}(a,c,\mu )f(z)\right) 
{\small +}\tfrac{p{\small -}\xi }{p}\left( z^{p}I_{\lambda ,\ell
}^{p,m}(a,c,\mu )f(z)\right) \prec \tfrac{1+Az}{1+Bz}{\small +}\tfrac{%
\lambda \xi }{\ell p}\tfrac{\left( A{\small -}B\right) z}{\left( 1{\small +}%
Bz\right) ^{2}}\text{. \ \ \ \ \ \ \ \ \ }  \tag{3.20}
\end{equation}

\noindent Then%
\begin{equation}
z^{p}I_{\lambda ,\ell }^{p,m}(a,c,\mu )f(z)\prec \tfrac{1+Az}{1+Bz}, 
\tag{3.21}
\end{equation}

\noindent and $\tfrac{1+Az}{1+Bz}$ is the best dominant of (3.21).\medskip

\noindent \textbf{Proof.} Uppon setting $q(z)=\dfrac{1+Az}{1+Bz},$ we see
that%
\begin{equation*}
1+\frac{zq^{\prime \prime }(z)}{q^{\prime }(z)}=\frac{1-Bz}{1+Bz},
\end{equation*}%
then, we get%
\begin{equation*}
\func{Re}\left\{ 1+\frac{zq^{\prime \prime }(z)}{q^{\prime }(z)}\right\} >%
\frac{1-\left\vert B\right\vert }{1+\left\vert B\right\vert }\text{ \ }%
\left( z\in U\right) .
\end{equation*}

\noindent Consequently, the hypotheses (3.15), (3.17) and (3.19) imply the
conditions (3.1), (3.3), and (3.5) respectively of Theorem 1. Therefore, the
assertion (3.21) follows from Theorem 1. The proof of Corollary 1 is
completed.\bigskip

\noindent Taking $p=A=1$ and $B=-1$ in Corollary 1,we obtain the following
corollary:

\noindent \textbf{Corollary 2.} Let $\xi \in 
%TCIMACRO{\U{2102} }%
%BeginExpansion
\mathbb{C}
%EndExpansion
^{\ast }$. Let the function $f\in \Sigma $. Suppose any one of the following
pairs of conditions is satisfied:%
\begin{equation}
\func{Re}\left\{ \dfrac{a-\mu }{\xi }\right\} >0,\text{ \ \ \ \ \ \ \ \ \ \
\ \ \ \ \ \ \ \ \ \ \ \ \ \ \ \ \ \ \ \ \ \ \ \ \ \ \ \ \ \ \ }  \tag{3.22}
\end{equation}%
\begin{equation}
\xi \left( zI_{\lambda ,\ell }^{m}(a{\small +}1,c,\mu )f(z)\right) {\small +}%
\left( 1{\small -}\xi \right) \left( zI_{\lambda ,\ell }^{m}(a,c,\mu
)f(z)\right) \prec \dfrac{1+z}{1-z}{\small +}\dfrac{\mu \xi }{a{\small -}\mu 
}\dfrac{2z}{\left( 1{\small -}z\right) ^{2}},  \tag{3.23}
\end{equation}

\noindent or

\begin{equation}
\func{Re}\left\{ \dfrac{c-\mu -1}{\xi }\right\} >0,\text{ \ \ \ \ \ \ \ \ \
\ \ \ \ \ \ \ \ \ \ \ \ \ \ \ \ \ \ \ \ \ \ \ \ \ \ \ \ \ \ }  \tag{3.24}
\end{equation}%
\begin{equation}
\xi \left( zI_{\lambda ,\ell }^{m}(a,c{\small -}1,\mu )f(z)\right) {\small +}%
\left( 1{\small -}\xi \right) \left( zI_{\lambda ,\ell }^{m}(a,c,\mu
)f(z)\right) \prec \dfrac{1+z}{1-z}{\small +}\dfrac{\mu \xi }{c{\small -}\mu 
{\small -}1}\dfrac{2z}{\left( 1{\small -}z\right) ^{2}},  \tag{3.25}
\end{equation}%
\noindent or%
\begin{equation}
\func{Re}\left\{ \dfrac{1}{\xi }\right\} >0,\text{ \ \ \ \ \ \ \ \ \ \ \ \ \
\ \ \ \ \ \ \ \ \ \ \ \ \ \ \ \ \ \ \ \ \ \ \ \ \ \ \ \ \ \ \ }  \tag{3.26}
\end{equation}%
\begin{equation}
\xi \left( zI_{\lambda ,\ell }^{m-1}(a,c,\mu )f(z)\right) {\small +}\left( 1%
{\small -}\xi \right) \left( zI_{\lambda ,\ell }^{m}(a,c,\mu )f(z)\right)
\prec \dfrac{1+z}{1-z}{\small +}\dfrac{\lambda \xi }{\ell }\dfrac{2z}{\left(
1{\small -}z\right) ^{2}}\text{. \ \ \ \ \ \ \ \ \ }  \tag{3.27}
\end{equation}

\noindent Then%
\begin{equation}
zI_{\lambda ,\ell }^{m}(a,c,\mu )f(z)\prec \dfrac{1+z}{1-z},  \tag{3.28}
\end{equation}

\noindent and $\tfrac{1+z}{1-z}$ is the best dominant of (3.28).\medskip

\noindent Taking $a=c$ and $m=0$ in Corollary 2,we obtain the following
corollary:

\noindent \textbf{Corollary 3. }Let $\xi \in 
%TCIMACRO{\U{2102} }%
%BeginExpansion
\mathbb{C}
%EndExpansion
^{\ast }$. Let the function $f\in \Sigma $. Suppose any one of the following
pairs of conditions is satisfied:%
\begin{equation}
\func{Re}\left\{ \dfrac{a-\mu }{\xi }\right\} >0,\text{ \ \ \ \ \ \ \ \ \ \
\ \ \ \ \ \ \ \ \ \ \ \ \ \ \ \ \ \ \ \ \ \ \ \ \ \ \ \ \ \ \ }  \tag{3.29}
\end{equation}%
\begin{equation}
\frac{\mu \xi }{a-\mu }z\left( zf(z)\right) ^{\prime }+zf(z)\prec \dfrac{1+z%
}{1-z}{\small +}\dfrac{\mu \xi }{a{\small -}\mu }\dfrac{2z}{\left( 1{\small -%
}z\right) ^{2}},  \tag{3.30}
\end{equation}

\noindent or

\begin{equation}
\func{Re}\left\{ \dfrac{c-\mu -1}{\xi }\right\} >0,\text{ \ \ \ \ \ \ \ \ \
\ \ \ \ \ \ \ \ \ \ \ \ \ \ \ \ \ \ \ \ \ \ \ \ \ \ \ \ \ \ }  \tag{3.31}
\end{equation}%
\begin{equation}
\frac{\mu \xi }{c-\mu -1}z\left( zf(z)\right) ^{\prime }+zf(z)\prec \dfrac{%
1+z}{1-z}{\small +}\dfrac{\mu \xi }{c{\small -}\mu {\small -}1}\dfrac{2z}{%
\left( 1{\small -}z\right) ^{2}},  \tag{3.32}
\end{equation}%
\noindent or%
\begin{equation}
\func{Re}\left\{ \dfrac{1}{\xi }\right\} >0,\text{ \ \ \ \ \ \ \ \ \ \ \ \ \
\ \ \ \ \ \ \ \ \ \ \ \ \ \ \ \ \ \ \ \ \ \ \ \ \ \ \ \ \ \ \ }  \tag{3.33}
\end{equation}%
\begin{equation}
\frac{\lambda \xi }{\ell }z\left( zf(z)\right) ^{\prime }+zf(z)\prec \dfrac{%
1+z}{1-z}{\small +}\dfrac{\lambda \xi }{\ell }\dfrac{2z}{\left( 1{\small -}%
z\right) ^{2}}\text{. \ \ \ \ \ \ \ \ \ }  \tag{3.34}
\end{equation}

\noindent Then%
\begin{equation}
zf(z)\prec \dfrac{1+z}{1-z},  \tag{3.35}
\end{equation}

\noindent and $\tfrac{1+z}{1-z}$ is the best dominant of (3.35).

\noindent Also, we introduce another subordination theorem as follows:

\noindent \textbf{Theorem 2.} Let $q(z)$ be a non zero univalent function in 
$U$ with $q(0)=1$. Let $\eta \in 
%TCIMACRO{\U{2102} }%
%BeginExpansion
\mathbb{C}
%EndExpansion
^{\ast }$ and $\tau ,\varkappa \in 
%TCIMACRO{\U{2102} }%
%BeginExpansion
\mathbb{C}
%EndExpansion
$ with $\tau +\varkappa \neq 0$. Let $f\in \Sigma _{p}$ and suppose that $f$
and $q$ satisfy the conditions:%
\begin{equation*}
\frac{\tau z^{p}I_{\lambda ,\ell }^{p,m}(a{\small +}1,c,\mu )f(z){\small %
+\varkappa }z^{p}I_{\lambda ,\ell }^{p,m}(a,c,\mu )f(z)}{\tau +\varkappa }%
\neq 0\text{ \ }\left( z\in U\right) ,
\end{equation*}

\noindent and%
\begin{equation}
\func{Re}\left\{ 1+\frac{zq^{\prime \prime }(z)}{q^{\prime }(z)}-\frac{%
zq^{\prime }(z)}{q(z)}\right\} >0\text{ \ }\left( z\in U\right) .  \tag{3.36}
\end{equation}

\noindent If%
\begin{equation}
\eta \left[ p+\frac{\tau z\left( I_{\lambda ,\ell }^{p,m}(a{\small +}1,c,\mu
)f(z)\overset{}{}\right) ^{\prime }{\small +\varkappa }z\left( I_{\lambda
,\ell }^{p,m}(a,c,\mu )f(z)\overset{}{}\right) ^{\prime }}{\tau I_{\lambda
,\ell }^{p,m}(a{\small +}1,c,\mu )f(z){\small +\varkappa }I_{\lambda ,\ell
}^{p,m}(a,c,\mu )f(z)}\right] \prec \frac{zq^{\prime }(z)}{q(z)},  \tag{3.37}
\end{equation}

\noindent then%
\begin{equation}
\left[ \frac{\tau z^{p}I_{\lambda ,\ell }^{p,m}(a{\small +}1,c,\mu )f(z)%
{\small +\varkappa }z^{p}I_{\lambda ,\ell }^{p,m}(a,c,\mu )f(z)}{\tau
+\varkappa }\right] ^{\eta }\prec q(z),  \tag{3.38}
\end{equation}

\noindent and $q(z)$ is the best dominant of (3.38).

\noindent \textbf{Proof.} In view of Lemma 2, we set%
\begin{equation*}
\theta (w)=0\text{ \ and\ \ }\varphi (w)=\frac{1}{w}\text{.}
\end{equation*}

\noindent thus%
\begin{equation*}
Q(z)=zq^{\prime }(z)\varphi (q(z))=\frac{zq^{\prime }(z)}{q(z)}\text{ \ and
\ }h(z)=Q(z).
\end{equation*}

\noindent By hypothesis (3.36), we note that $Q(z)$ is univalent, moreover%
\begin{equation*}
\func{Re}\left\{ \frac{zQ^{\prime }(z)}{Q(z)}\right\} =\func{Re}\left\{ 
\frac{z\left( \frac{zq^{\prime }(z)}{q(z)}\right) ^{\prime }}{\frac{%
zq^{\prime }(z)}{q(z)}}\right\} =\func{Re}\left\{ 1{\small +}\frac{%
zq^{\prime \prime }(z)}{q^{\prime }(z)}{\small -}\frac{zq^{\prime }(z)}{q(z)}%
\right\} >0\text{ \ }\left( z\in U\right) ,
\end{equation*}%
then function $Q(z)$ is also starlike in $U$. We furthermore get that%
\begin{equation*}
\func{Re}\left\{ \frac{zh^{\prime }(z)}{Q(z)}\right\} >0\text{ \ }\left(
z\in U\right) .
\end{equation*}

\noindent Next, let the function $p$ be defined by%
\begin{equation}
p(z)=\left[ \frac{\tau z^{p}I_{\lambda ,\ell }^{p,m}(a{\small +}1,c,\mu )f(z)%
{\small +\varkappa }z^{p}I_{\lambda ,\ell }^{p,m}(a,c,\mu )f(z)}{\tau
+\varkappa }\right] ^{\eta }\text{ \ }\left( z\in U\right) .  \tag{3.39}
\end{equation}

\noindent Then $p$ is analytic in $U$, $p(0)=q(0)=1$ and%
\begin{equation}
\frac{zp^{\prime }(z)}{p(z)}=\eta \left[ p+\frac{\tau z\left( I_{\lambda
,\ell }^{p,m}(a{\small +}1,c,\mu )f(z)\overset{}{}\right) ^{\prime }{\small %
+\varkappa }z\left( I_{\lambda ,\ell }^{p,m}(a,c,\mu )f(z)\overset{}{}%
\right) ^{\prime }}{\tau I_{\lambda ,\ell }^{p,m}(a{\small +}1,c,\mu )f(z)%
{\small +\varkappa }I_{\lambda ,\ell }^{p,m}(a,c,\mu )f(z)}\right] . 
\tag{3.40}
\end{equation}

\noindent Using (3.40) in (3.37), we have%
\begin{equation*}
\frac{zp^{\prime }(z)}{p(z)}\prec \frac{zq^{\prime }(z)}{q(z)},
\end{equation*}%
\noindent which is also equivalent to%
\begin{equation*}
zp^{\prime }(z)\varphi (p(z))\prec zq^{\prime }(z)\varphi (q(z)),
\end{equation*}%
\noindent or%
\begin{equation*}
\theta (p(z))+zp^{\prime }(z)\varphi (p(z))\prec \theta (q(z))+zq^{\prime
}(z)\varphi (q(z)).
\end{equation*}%
\noindent Therefore, by Lemma 2, we have%
\begin{equation*}
p(z)\prec q(z),
\end{equation*}

\noindent and $q(z)$ is the best dominant. This is precisely the assertion
in (3.38). The proof of Theorem 2 is completed.\bigskip

\noindent Taking $\tau =0$, $\varkappa =1$ and $q(z)=\tfrac{1+Az}{1+Bz}$ in
Theorem 2, we obtain the following corollary.

\noindent \textbf{Corollary 4. }Let $\eta \in 
%TCIMACRO{\U{2102} }%
%BeginExpansion
\mathbb{C}
%EndExpansion
^{\ast }$. Let $f\in \Sigma _{p}$ and suppose that $f$ satisfies the
conditions:%
\begin{equation*}
z^{p}I_{\lambda ,\ell }^{p,m}(a,c,\mu )f(z)\neq 0\text{ \ }\left( z\in
U\right) ,
\end{equation*}%
\noindent if%
\begin{equation}
\eta \left[ p+\frac{z\left( I_{\lambda ,\ell }^{p,m}(a,c,\mu )f(z)\overset{}{%
}\right) ^{\prime }}{I_{\lambda ,\ell }^{p,m}(a,c,\mu )f(z)}\right] \prec 
\frac{\left( A-B\right) z}{\left( 1+Az\right) \left( 1+Bz\right) }, 
\tag{3.41}
\end{equation}

\noindent then%
\begin{equation}
\left[ z^{p}I_{\lambda ,\ell }^{p,m}(a,c,\mu )f(z)\right] ^{\eta }\prec 
\dfrac{1+Az}{1+Bz},  \tag{3.42}
\end{equation}

\noindent and $\tfrac{1+Az}{1+Bz}$ is the best dominant of (3.42).\bigskip

\noindent Taking $p=A=1$ and $B=-1$ in Corollary 4, we obtain the following
corollary:

\noindent \textbf{Corollary 5. }Let $\eta \in 
%TCIMACRO{\U{2102} }%
%BeginExpansion
\mathbb{C}
%EndExpansion
^{\ast }$. Let $f\in \Sigma $ and suppose that $f$ satisfies the conditions:%
\begin{equation*}
zI_{\lambda ,\ell }^{m}(a,c,\mu )f(z)\neq 0\text{ \ }\left( z\in U\right) ,
\end{equation*}%
\noindent if%
\begin{equation}
\eta \left[ 1+\frac{z\left( I_{\lambda ,\ell }^{m}(a,c,\mu )f(z)\overset{}{}%
\right) ^{\prime }}{I_{\lambda ,\ell }^{m}(a,c,\mu )f(z)}\right] \prec \frac{%
2z}{\left( 1-z^{2}\right) },  \tag{3.43}
\end{equation}

\noindent then%
\begin{equation}
\left[ zI_{\lambda ,\ell }^{m}(a,c,\mu )f(z)\right] ^{\eta }\prec \dfrac{1+z%
}{1-z},  \tag{3.44}
\end{equation}

\noindent and $\tfrac{1+z}{1-z}$ is the best dominant of (3.44).\bigskip

\noindent Taking $a=c$, $\eta =1$ and $m=0$ in Corollary 5, we obtain the
following corollary:

\noindent \textbf{Corollary 6. }Let $f\in \Sigma $ and suppose that $f$
satisfies the conditions:%
\begin{equation*}
zf(z)\neq 0\text{ \ }\left( z\in U\right) ,
\end{equation*}%
\noindent if%
\begin{equation}
1+\frac{zf^{\prime }(z)}{f(z)}\prec \frac{2z}{\left( 1-z^{2}\right) }, 
\tag{3.45}
\end{equation}

\noindent then%
\begin{equation}
zf(z)\prec \dfrac{1+z}{1-z},  \tag{3.46}
\end{equation}

\noindent and $\dfrac{1+z}{1-z}$ is the best dominant of (3.46).\bigskip

\noindent Taking $\tau =1$, $\varkappa =0$ and $q(z)=\dfrac{1+Az}{1+Bz}$ in
Theorem 2, we obtain the following corollary.

\noindent \textbf{Corollary 7. }Let $\eta \in 
%TCIMACRO{\U{2102} }%
%BeginExpansion
\mathbb{C}
%EndExpansion
^{\ast }$. Let $f\in \Sigma _{p}$ and suppose that $f$ satisfies the
conditions:%
\begin{equation*}
z^{p}I_{\lambda ,\ell }^{p,m}(a{\small +}1,c,\mu )f(z)\neq 0\text{ \ }\left(
z\in U\right) ,
\end{equation*}%
\noindent if%
\begin{equation}
\eta \left[ p+\frac{z\left( I_{\lambda ,\ell }^{p,m}(a+1,c,\mu )f(z)\overset{%
}{}\right) ^{\prime }}{I_{\lambda ,\ell }^{p,m}(a+1,c,\mu )f(z)}\right]
\prec \frac{\left( A-B\right) z}{\left( 1+Az\right) \left( 1+Bz\right) }, 
\tag{3.47}
\end{equation}

\noindent then%
\begin{equation}
\left[ z^{p}I_{\lambda ,\ell }^{p,m}(a+1,c,\mu )f(z)\right] ^{\eta }\prec 
\dfrac{1+Az}{1+Bz},  \tag{3.48}
\end{equation}

\noindent and $\tfrac{1+Az}{1+Bz}$ is the best dominant of (3.48).\bigskip

\noindent Taking $A=p=1$ and $B=-1$ in Corollary 7, we obtain the following
corollary.

\noindent \textbf{Corollary 8. }Let $\eta \in 
%TCIMACRO{\U{2102} }%
%BeginExpansion
\mathbb{C}
%EndExpansion
^{\ast }$. Let $f\in \Sigma $ and suppose that $f$ satisfies the conditions:%
\begin{equation*}
zI_{\lambda ,\ell }^{m}(a{\small +}1,c,\mu )f(z)\neq 0\text{ \ }\left( z\in
U\right) ,
\end{equation*}%
\noindent if%
\begin{equation}
\eta \left[ 1+\frac{z\left( I_{\lambda ,\ell }^{m}(a+1,c,\mu )f(z)\overset{}{%
}\right) ^{\prime }}{I_{\lambda ,\ell }^{m}(a+1,c,\mu )f(z)}\right] \prec 
\frac{2z}{\left( 1-z^{2}\right) },  \tag{3.49}
\end{equation}

\noindent then%
\begin{equation}
\left[ zI_{\lambda ,\ell }^{m}(a+1,c,\mu )f(z)\right] ^{\eta }\prec \dfrac{%
1+z}{1-z},  \tag{3.50}
\end{equation}

\noindent and $\tfrac{1+z}{1-z}$ is the best dominant of (3.50).\bigskip

\noindent Taking $a=c$, $\eta =1$ and $m=0$ in Corollary 8, we obtain the
following corollary:

\noindent \textbf{Corollary 9. }Let $f\in \Sigma $ and suppose that $f$
satisfies the conditions:%
\begin{equation*}
z^{2}f^{\prime }(z)+\frac{a}{\mu }zf(z)\neq 0\text{ \ }\left( z\in U\right) ,
\end{equation*}%
\noindent if%
\begin{equation}
1+\frac{z\left( z^{2}f^{\prime }(z)+\frac{a}{\mu }zf(z)\overset{}{}\right)
^{\prime }}{z^{2}f^{\prime }(z)+\frac{a}{\mu }zf(z)}\prec \frac{2z}{\left(
1-z^{2}\right) },  \tag{3.51}
\end{equation}

\noindent then%
\begin{equation}
\frac{\mu }{a-\mu }\left( z^{2}f^{\prime }(z)+\frac{a}{\mu }zf(z)\right)
\prec \dfrac{1+z}{1-z},  \tag{3.52}
\end{equation}

\noindent and $\tfrac{1+z}{1-z}$ is the best dominant of (3.52).\bigskip

\noindent Another theorem is introduced as follows:

\noindent \textbf{Theorem 3.} Let $\eta \in 
%TCIMACRO{\U{2102} }%
%BeginExpansion
\mathbb{C}
%EndExpansion
^{\ast }$ and $\zeta ,\tau ,\varkappa \in 
%TCIMACRO{\U{2102} }%
%BeginExpansion
\mathbb{C}
%EndExpansion
$ with $\tau +\varkappa \neq 0$. Let $q(z)$ be a univalent function in $U$
with $q(0)=1$ and%
\begin{equation}
\func{Re}\left\{ 1+\frac{zq^{\prime \prime }(z)}{q^{\prime }(z)}\right\}
>\max \left\{ 0,-\func{Re}\left\{ \zeta \right\} \right\} \text{ \ }\left(
z\in U\right) .  \tag{3.53}
\end{equation}%
Let $f\in \Sigma _{p}$ and suppose that $f$ satisfies the condition%
\begin{equation*}
\frac{\tau z^{p}I_{\lambda ,\ell }^{p,m}(a{\small +}1,c,\mu )f(z){\small %
+\varkappa }z^{p}I_{\lambda ,\ell }^{p,m}(a,c,\mu )f(z)}{\tau +\varkappa }%
\neq 0\text{ \ }\left( z\in U\right) ,
\end{equation*}

\noindent Set%
\begin{eqnarray}
\Omega (z) &=&\left[ \tfrac{\tau z^{p}I_{\lambda ,\ell }^{p,m}(a{\small +}%
1,c,\mu )f(z){\small +\varkappa }z^{p}I_{\lambda ,\ell }^{p,m}(a,c,\mu )f(z)%
}{\tau +\varkappa }\right] ^{\eta }  \notag \\
&&\cdot \left[ \zeta +\eta \left( \tfrac{\tau z\left( I_{\lambda ,\ell
}^{p,m}(a{\small +}1,c,\mu )f(z)\right) ^{\prime }{\small +\varkappa }%
z\left( I_{\lambda ,\ell }^{p,m}(a,c,\mu )f(z)\right) ^{\prime }}{\tau
I_{\lambda ,\ell }^{p,m}(a{\small +}1,c,\mu )f(z){\small +\varkappa }%
I_{\lambda ,\ell }^{p,m}(a,c,\mu )f(z)}+p\right) \right] .  \TCItag{3.54}
\end{eqnarray}

\noindent If%
\begin{equation}
\Omega (z)\prec \zeta q(z)+zq^{\prime }(z),  \tag{3.55}
\end{equation}

\noindent then%
\begin{equation}
\left[ \frac{\tau z^{p}I_{\lambda ,\ell }^{p,m}(a{\small +}1,c,\mu )f(z)%
{\small +\varkappa }z^{p}I_{\lambda ,\ell }^{p,m}(a,c,\mu )f(z)}{\tau
+\varkappa }\right] ^{\eta }\prec q(z),  \tag{3.56}
\end{equation}

\noindent and $q(z)$ is the best dominant of (3.56).

\noindent Proof. In view of Lemma 2, we set%
\begin{equation*}
\theta (w)=\zeta w\text{ \ and\ \ }\varphi (w)=1\text{ \ }(w\in 
%TCIMACRO{\U{2102} }%
%BeginExpansion
\mathbb{C}
%EndExpansion
)\text{,}
\end{equation*}

\noindent thus%
\begin{equation*}
Q(z)=zq^{\prime }(z)\varphi (q(z))=zq^{\prime }(z)\text{ \ and \ }h(z)=\zeta
q(z)+zq^{\prime }(z).
\end{equation*}

\noindent Then, we note that $Q(z)$ is univalent. Moreover, using (3.53), we
find that%
\begin{equation*}
\func{Re}\left\{ \frac{zQ^{\prime }(z)}{Q(z)}\right\} =\func{Re}\left\{ 
\frac{z\left( zq^{\prime }(z)\right) ^{\prime }}{zq^{\prime }(z)}\right\} =%
\func{Re}\left\{ 1{\small +}\frac{zq^{\prime \prime }(z)}{q^{\prime }(z)}%
\right\} >0\text{ \ }\left( z\in U\right) ,
\end{equation*}%
then function $Q(z)$ is also starlike in $U$. Also, using (3.53), we get that%
\begin{equation*}
\func{Re}\left\{ \frac{zh^{\prime }(z)}{Q(z)}\right\} =\func{Re}\left\{
1+\zeta +\frac{zq^{\prime \prime }(z)}{q^{\prime }(z)}\right\} >0\text{ \ }%
\left( z\in U\right) .
\end{equation*}%
\noindent Furthermore, by using the expression of $p(z)$ defined by (3.39)
and the expression of $zp\prime (z)$ defined by (3.40) we have%
\begin{eqnarray*}
\theta (p(z)){\small +}zp^{\prime }(z)\varphi (p(z)) &=&\zeta
p(z)+zp^{\prime }(z) \\
&=&\left[ \tfrac{\tau z^{p}I_{\lambda ,\ell }^{p,m}(a{\small +}1,c,\mu )f(z)%
{\small +\varkappa }z^{p}I_{\lambda ,\ell }^{p,m}(a,c,\mu )f(z)}{\tau
+\varkappa }\right] ^{\eta } \\
&&\cdot \left[ \zeta {\small +}\eta \left( \tfrac{\tau z\left( I_{\lambda
,\ell }^{p,m}(a{\small +}1,c,\mu )f(z)\right) ^{\prime }{\small +\varkappa }%
z\left( I_{\lambda ,\ell }^{p,m}(a,c,\mu )f(z)\right) ^{\prime }}{\tau
I_{\lambda ,\ell }^{p,m}(a{\small +}1,c,\mu )f(z){\small +\varkappa }%
I_{\lambda ,\ell }^{p,m}(a,c,\mu )f(z)}{\small +}p\right) \right] \\
&=&\Omega (z).
\end{eqnarray*}%
\noindent The hypothesis (3.55) is now equivalent to%
\begin{equation*}
\zeta p(z)+zp^{\prime }(z)\prec \zeta q(z)+zq^{\prime }(z),
\end{equation*}%
\noindent or%
\begin{equation*}
\theta (p(z)){\small +}zp^{\prime }(z)\varphi (p(z))\prec \theta (q(z))%
{\small +}zq^{\prime }(z)\varphi (q(z)).
\end{equation*}

\noindent Finally, an application of Lemma 2 yields%
\begin{equation*}
p(z)\prec q(z)
\end{equation*}

\noindent and $q(z)$ is the best dominant. This is precisely the assertion
in (3.56). The proof of Theorem 3 is completed.\bigskip

\noindent Taking $\tau {\small =}0$, $\varkappa {\small =}1$ and $q(z)%
{\small =}\dfrac{1{\small +}Az}{1{\small +}Bz}$ in Theorem 3, we obtain the
following corollary.

\noindent \textbf{Corollary 10. }Let $\eta \in 
%TCIMACRO{\U{2102} }%
%BeginExpansion
\mathbb{C}
%EndExpansion
^{\ast }$ and $\zeta =\frac{\left\vert B\right\vert -1}{\left\vert
B\right\vert +1}$. Let $f\in \Sigma _{p}$ and suppose that $f$ satisfies the
conditions%
\begin{equation*}
z^{p}I_{\lambda ,\ell }^{p,m}(a,c,\mu )f(z)\neq 0\text{ \ }\left( z\in
U\right) ,
\end{equation*}

\noindent and%
\begin{equation}
\left[ z^{p}I_{\lambda ,\ell }^{p,m}(a,c,\mu )f(z)\underset{}{}\right]
^{\eta }\cdot \left[ \zeta {\small +}\eta \left( p{\small +}\tfrac{z\left(
I_{\lambda ,\ell }^{p,m}(a,c,\mu )f(z)\underset{}{}\right) ^{\prime }}{%
I_{\lambda ,\ell }^{p,m}(a,c,\mu )f(z)}\right) \right] \prec \zeta \dfrac{1%
{\small +}Az}{1{\small +}Bz}{\small +}\frac{(A{\small -}B)z}{(1{\small +}%
Bz)^{2}},  \tag{3.57}
\end{equation}

\noindent then%
\begin{equation}
\left[ z^{p}I_{\lambda ,\ell }^{p,m}(a,c,\mu )f(z)\underset{}{}\right]
^{\eta }\prec \dfrac{1{\small +}Az}{1{\small +}Bz},  \tag{3.58}
\end{equation}

\noindent and $\dfrac{1{\small +}Az}{1{\small +}Bz}$ is the best dominant of
(3.58).\bigskip

\noindent Taking $p=A=1$, $B=-1$ and $a=c$ in Corollary 10, we obtain the
following corollary.

\noindent \textbf{Corollary 11. }Let $\eta \in 
%TCIMACRO{\U{2102} }%
%BeginExpansion
\mathbb{C}
%EndExpansion
^{\ast }$. Let $f\in \Sigma $ and suppose that $f$ satisfies the conditions%
\begin{equation*}
zf(z)\neq 0\text{ \ }\left( z\in U\right) ,
\end{equation*}

\noindent and%
\begin{equation}
\left[ zf(z)\underset{}{}\right] ^{\eta }\cdot \left[ \eta \left( 1+\dfrac{%
zf^{\prime }(z)}{f(z)}\right) \right] \prec \frac{2z}{(1-z)^{2}},  \tag{3.59}
\end{equation}

\noindent then%
\begin{equation}
\left[ zf(z)\underset{}{}\right] ^{\eta }\prec \dfrac{1+z}{1-z},  \tag{3.60}
\end{equation}

\noindent and $\dfrac{1{\small +}z}{1{\small -}z}$ is the best dominant of
(3.60).

\noindent \textbf{Remark 1. }The result obtained in Corollary 11 coincides
with the recent result due to Mishra et al. [14, Corollary 4.9].\bigskip

\noindent Taking $\eta =1$ in Corollary 11, we obtain the following
corollary.

\noindent \textbf{Corollary 12. }Let $f\in \Sigma $ and suppose that $f$
satisfies the conditions%
\begin{equation}
zf(z)\neq 0\text{ \ }\left( z\in U\right) ,  \tag{3.61}
\end{equation}%
\noindent and%
\begin{equation}
zf(z)+z^{2}f^{\prime }(z)\prec \frac{2z}{(1-z)^{2}},  \tag{3.62}
\end{equation}%
\noindent then%
\begin{equation}
zf(z)\prec \dfrac{1+z}{1-z},  \tag{3.63}
\end{equation}

\noindent and $\dfrac{1{\small +}z}{1{\small -}z}$ is the best dominant of
(3.63).

\noindent Taking $\tau {\small =}1$, $\varkappa {\small =}0$ and $q(z)%
{\small =}\tfrac{1{\small +}Az}{1{\small +}Bz}$ in Theorem 3, we obtain the
following corollary.

\noindent \textbf{Corollary 13. }Let $\eta \in 
%TCIMACRO{\U{2102} }%
%BeginExpansion
\mathbb{C}
%EndExpansion
^{\ast }$ and $\zeta =\frac{\left\vert B\right\vert -1}{\left\vert
B\right\vert +1}$. Let $f\in \Sigma _{p}$ and suppose that $f$ satisfies the
conditions%
\begin{equation}
z^{p}I_{\lambda ,\ell }^{p,m}(a{\small +}1,c,\mu )f(z)\neq 0\text{ \ }\left(
z\in U\right) ,  \tag{3.64}
\end{equation}

\noindent and%
\begin{equation}
\left[ z^{p}I_{\lambda ,\ell }^{p,m}(a{\small +}1,c,\mu )f(z)\underset{}{}%
\right] ^{\eta }\cdot \left[ \zeta +\eta \left( \tfrac{z\left( I_{\lambda
,\ell }^{p,m}(a{\small +}1,c,\mu )f(z)\right) ^{\prime }}{I_{\lambda ,\ell
}^{p,m}(a{\small +}1,c,\mu )f(z)}+p\right) \right] \prec \zeta \tfrac{1%
{\small +}Az}{1{\small +}Bz}{\small +}\tfrac{(A{\small -}B)z}{(1{\small +}%
Bz)^{2}},  \tag{3.65}
\end{equation}

\noindent then%
\begin{equation}
\left[ z^{p}I_{\lambda ,\ell }^{p,m}(a{\small +}1,c,\mu )f(z)\right] ^{\eta
}\prec \dfrac{1{\small +}Az}{1{\small +}Bz},  \tag{3.66}
\end{equation}

\noindent and $q(z)$ is the best dominant of (3.66).\bigskip

\noindent Taking $p=A=\eta =1$, $B=-1$ and $a=c$ in Corollary 13, we obtain
the following corollary.

\noindent \textbf{Corollary 14. }Let $f\in \Sigma $ and suppose that $f$
satisfies the conditions

\begin{equation}
z^{2}f^{\prime }(z)+\frac{a}{\mu }zf(z)\neq 0\text{ \ }\left( z\in U\right) ,
\tag{3.67}
\end{equation}

\noindent and%
\begin{equation}
\frac{\mu z}{a-\mu }\left( z\left[ zf^{\prime }(z)+\dfrac{a}{\mu }f(z)\right]
\underset{}{}\right) ^{\prime }\prec \frac{2z}{(1{\small -}z)^{2}} 
\tag{3.68}
\end{equation}

\noindent then%
\begin{equation}
\frac{\mu }{a-\mu }\left( z^{2}f^{\prime }(z)+\frac{a}{\mu }zf(z)\right)
\prec \dfrac{1+z}{1-z},  \tag{3.69}
\end{equation}

\noindent and $\tfrac{1+z}{1-z}$ is the best dominant of (3.69).

\noindent \textbf{Remark 2. }Specializing the parameters in Theorems 1,2 and
3 as mentioned before, we can obtain the corresponding subordination
properties of Liu-Srivastava operator [11], Cho-Kwon-Srivastava operator
[5], Yuan-Liu-Srivastava operator [20], Uralegaddi-Somanatha operator [19],
and others.

\end{document}